\documentclass[12 pt]{article}
\newtheorem{theorem}{Theorem}
\begin{document}
\title{A limit process for a sequence of partial sums of residuals of a simple regression against order statistics
with Markov-modulated noise}
\author{Artyom Kovalevskii \thanks{E-mail: pandorra@ngs.ru},  Evgeny Shatalin \thanks{E-mail: sh\_e\_v\_89@list.ru}}
\date{}
\maketitle

\begin{abstract}
We consider a simple regression model where a regressor is composed
of order statistics and a noise is Markov-modulated. We introduce an
empirical bridge of regression residuals and prove its weak
convergence to a centered Gaussian process.
\end{abstract}


\section{Introduction and main results}
Brown et al. (1975) proposed a test for change of regression at
unknown time. Their approach is based on computation of recursive
residuals. MacNeill (1978) studied a linear regression against
values of continuously differentiable functions. He obtained limit
processes for sequences of partial sums of regression residuals.
Later Bischoff (1997) showed that the MacNeill's theorem holds in
more general setting, namely for continuous regressor functions. Aue
et al. (2008) introduced a new test for polynomial regression
functions which is analogous to the classical likelihood test. Stute
(1997) proposed a class of tests that are based on regression
residuals. His general approach also allows to analyse models with order
statistics regressors.

We consider another model of a simple linear regression against
order statistics where the noise is Markov-modulated, and analyse a
limit process for sums of regression residuals.

To define the model, we introduce 3 mutually independent families of
random variables:

1) $\{\varepsilon_i^v, i\geq 1, 1\leq v\leq M\}$, a family of
independent random variables where $\{\varepsilon_i^v, i\geq 1\}$
are identically distributed for each $v$, ${\bf
E}\varepsilon_1^v=0$, ${\bf D}\varepsilon_1^v=\sigma_v^2\geq0$ and
$\sum_{v=1}^{M}\sigma_v^2>0$;

2) $\{\xi_i\}_{i=1}^{\infty}$, a sequence of i.i.d. random variables
with distribution function $F$ and finite positive variance ${\bf
Var}\xi$;

3) $\{V_i\}_{i=1}^{\infty}$, an irreducible aperiodic Markov chain
on the state space $\{1,\dots,M\}$ with stationary distribution
$\{\pi_i\}_{i=1}^{M}$.

For any $n=1,2,\dots$, let $X_{ni}=\xi_{i:n}$ be the $i$-th order
statistic of the first $n$ random variables $\xi_1,\dots,\xi_n$,
where, in particular, $X_{n1}=\min_{1\leq i\leq n}\xi_i$ and
$X_{nn}=\max_{1\leq i\leq n}\xi_i$.

In this article, we consider the following regression model:

\[
Y_{ni}=a+bX_{ni}+\varepsilon_i^{V_i},~n\geq 1, i=1,\dots,n.
\]

For this model, we introduce an {\it empirical bridge} and show its weak convergence
to a centered Gaussian process.

The novelty of our model lies in consideration both ordered
regressors and Markov-modulated noise.

Let
\[
\widehat{b}_n=\frac{\overline{XY}-\overline{X}~\overline{Y}}{\overline{X^2}-\overline{X}^2},
~\widehat{a}_n=\overline{Y}-\widehat{b}_n\,\overline{X}.
\]
be the classical Gauss-Markov estimators for $a$ and $b$.

Define {\it fitted values} $\{\widehat{Y}_{ni} \}$, {\it regression
residuals} $\{\widehat{\varepsilon}_{ni}\}$ and their {\it partial
sums} $\{\widehat{\Delta}^0_{ni}\}$, by $\widehat{Y}_{ni} =
\widehat{a}_n+\widehat{b}_n X_{ni}$,
$\widehat{\varepsilon}_{ni}=Y_{ni} - \widehat{Y}_{ni}$ and
$\widehat{\Delta}_{ni}^0=~\widehat{\varepsilon}_{n1}+~\ldots+~\widehat{\varepsilon}_{ni}$,
for $1\leq i \leq n, \ n\geq 1$.

In what follows, we write for short: $Y_{ni}=Y_i$, $X_{ni}=X_i$,
$\widehat{Y}_{ni}=\widehat{Y}_{i}$,
$\widehat{\varepsilon}_{ni}=\widehat{\varepsilon}_{i}$ and
$\widehat{\Delta}^0_{ni}=\widehat{\Delta}^0_{i}$.

A {\it random polygon} $Z_n$ is a piecewise linear function with
nodes $(k/n, \widehat{\Delta}_k^0/\sigma\sqrt{n})$, for
$k=0,\dots,n$.

Further, an {\it empirical bridge} is a random
polygon $\widehat{Z}_n$ with nodes
$(k/n,(\widehat{\Delta}_k^0-k/n\widehat{\Delta}_n^0)/\sqrt{n\widehat{\sigma^2}})$
where
$\widehat{\sigma^2}=\overline{\widehat{\varepsilon}^2}-(\overline{\widehat{\varepsilon}})^2$
is an estimator of variance
$\sigma^2=\sum_{v=1}^{M}\sigma_v^2\pi_v$.

Let $GL_F(t)=\int\limits_{0}^{t}F^{-1}(s)\,ds$  be the {\it
theoretical general Lorenz curve} (Gastwirth, 1971; Davydov and
Zitikis, 2004) where $F^{-1}(s)=\sup\{x:~F(x)<s\}$ is the inverse of
distribution function $F(x)$. Let $GL_F^0(t)=GL_F(t)-tGL_F(1)$ be
its centered version. Similarly, let
$GL_n(t)=\frac{1}{n}\sum_{i=1}^{[nt]}\xi_{i:n}$  be the {\it
empirical Lorenz curve}. Goldie  (1977) showed that, as $n \to
\infty$, the empirical Lorenz curve converges a.s. to the
theoretical curve in the uniform metric, i.e. $\sup_{t \in {\bf R}}
| GL_n(t) - GL_F(t) | \to 0$ a.s.

Now we formulate the main result of the paper.

\begin{theorem} Both the random polygon
${Z_n}$ and the empirical bridge $\widehat{Z}_n$ converge weakly, as
$n\to\infty$, to the centered Gaussian process $Z_F^0$ with
covariance kernel, $K_F^0(t,s)$, given by
\[
K_F^0(t,s)=\min\{t,s\}-ts-\frac{GL_F^0(t)GL_F^0(s)}{{\bf Var}
\xi_1},~t, s\in[0,1].
\]
Here weak convergence holds in the space $C(0,1)$ of continuous
functions on [0,1] endowed by the uniform metric.
\end{theorem}

When the Markov chain degenerates, our model is a very particular
case of Stute (1997). Kovalevskii (2013) used this particular model
to analyse dependence of a car price on a production year.

In what follows, notation $\stackrel{\bf p}{\to}$ states for
convergence in probability.

\section{Proof of Theorem 1}
\

Let $X_i^0=X_i-\overline{X}$, $\varepsilon_i^0=\varepsilon_i^{V_i} -
\overline{\varepsilon}$ where
$\overline{\varepsilon}=\sum_{i=1}^{n}\varepsilon_i^{V_i}$.

The proof includes five steps. In the first step, we show that, in
the formulae under consideration, the sum
$\sum\limits_{i=1}^{n}\frac{\varepsilon_i^0X_i^0}{\sqrt{n}}$ may be
replaced by the sum $\sum\limits_{i=1}^{n}\frac{\varepsilon_i^0{\bf
E}X_i^0}{\sqrt{n}}$. Secondly, we prove weak convergence of a
normalized vector with coordinates
$(\widehat{\Delta}_{k_1}^0,\dots,\widehat{\Delta}_{k_m}^0)$ to a
normalized vector with coordinates
$(\Delta_{k_1}^0,\dots,\Delta_{k_m}^0)$ where $\Delta_{k_i}^0$ are
defined below. Then we prove weak convergence of finite-dimensional
distributions. The fourth step contains a proof of relative
compactness of the family $\{Z_n(t), 0\leq~t \leq 1\}$. We complete
with a proof of convergence of sample variance $\widehat{\sigma^2}$
to variance~$\sigma^2$.

{\bf Step 1}

Note that
\begin{equation}\label{1}
\widehat{\Delta}_k^0=\sum_{i=1}^{k} \left(\varepsilon_i^0
 - \frac{\overline{X^0 \varepsilon^0}}{\overline{(X^0)^2}}X_i^0\right).
\end{equation}

We show that
\begin{equation}\label{2}
\frac{1}{\sqrt{n}}\left(\sum\limits_{i=1}^{n}\varepsilon_i^0X_i^0
-\sum\limits_{i=1}^{n}\varepsilon_i^0{\bf
E}X_i^0\right)\stackrel{{\bf p}}{\to}0.
\end{equation}

Indeed,
\[
{\bf
P}\left\{\left|\frac{1}{\sqrt{n}}\sum\limits_{i=1}^{n}\varepsilon_i^0(X_i^0-{\bf
E}X_i^0)\right|\geq\delta\right\}\leq\frac{{\bf
Var}\sum\limits_{i=1}^{n}\varepsilon_i^0(X_i^0-{\bf
E}X_i^0)}{n\delta^2}
\]
\[
=\frac{\sum\limits_{i=1}^{n}{\bf Var}\varepsilon_i^{V_i}{\bf
Var}X_i^0}{n\delta^2}-\frac{2\sum\limits_{i,j=1}^{n}{\bf
Var}\varepsilon_i^{V_i}{\bf
cov}(X_i^0,X_j^0)}{n^2\delta^2}+\frac{\sum\limits_{k=1}^{n}{\bf
Var}\varepsilon_k^{V_k}\sum\limits_{i,j=1}^{n}{\bf
cov}(X_i^0,X_j^0)}{n^3\delta^2}.
\]
The last equality is correct because of the following equalities
\[
{\bf E}\varepsilon_i^{V_i}\varepsilon_j^{V_j}={\bf E}({\bf
E}\varepsilon_i^{v_i}\varepsilon_j^{v_j}|V_i=v_i,V_j=v_j) = \left\{
\begin{array}{ll}
{\bf
Var}\varepsilon_i^{V_i}, & i= j;\\
0, & i \neq j,
\end{array}
\right.
\]
\[
{\bf Var}\sum_{i=1}^{n}X_i\varepsilon_i^{V_i}={\bf E}({\bf
Var}\sum_{i=1}^{n}X_i\varepsilon_i^{v_i}|V_i=v_i,i=1,\dots,n)
\]
\[
={\bf E}(\sum_{i=1}^{n}{\bf Var}X_i{\bf
Var}\varepsilon_i^{v_i}|V_i=v_i,i=1,\dots,n)=\sum_{i=1}^{n}{\bf
Var}X_i{\bf Var}\varepsilon_i^{V_i},
\]
as $\{\varepsilon_i^{v_i}\}$ are i.i.d and do not depend on
$\{X_i\}$.

Theorem 1 (H\"{o}effding, 1953) implies
$\frac{1}{n}\sum\limits_{i=1}^{n}{\bf Var}X_i\to0$ as $n\to\infty$.

Note that 
 ${\bf Var}\overline{X}={\bf Var}\xi_1/n$,
 $\frac{1}{n}\sum_{i,j=1}^{n}{\bf cov}(X_i,X_j)=\frac{1}{n}{\bf
Var}\sum_{i=1}^{n}X_i={\bf Var}\xi_1$.

Prove that 
$\frac{1}{n}\sum\limits_{i=1}^{n}2|{\bf
cov}(X_i,\overline{X})|\to 0$
and
$\frac{1}{n}\sum\limits_{i=1}^{n}{\bf Var}X_i^0\to0$.
 
The sum of covariances admits the follows upper bound.
\[
\frac{1}{n}\sum\limits_{i=1}^{n}2|{\bf
cov}(X_i,\overline{X})|\leq\frac{1}{n}\sum\limits_{i=1}^{n}2\sqrt{{\bf
Var}X_i{\bf
Var}\overline{X}}\leq
\]
\[
\frac{1}{n}\sum\limits_{i=1}^{n}2\left(\frac{1+{\bf
Var}X_i}{2}\right)\sqrt{\frac{{\bf Var}X_i}{n}}
=\frac{1}{n\sqrt{n}}\sum\limits_{i=1}^{n}\Big(\sqrt{{\bf
Var}X_i}(1+{\bf Var}X_i)\Big)\leq
\]
\[
n^{-3/2}\sum_{i=1}^n \frac{1+{\bf Var} X_i}{2} + 
n^{-3/2}\left(\sum_{i=1}^n {\bf Var} X_i \right)^{3/2} \to 0.
\].

\[
\frac{1}{n}\sum\limits_{i=1}^{n}{\bf Var}X_i^0 =
\frac{1}{n}\sum_{i=1}^{n}\Big({\bf
Var}X_i+{\bf Var}\overline{X}-2{\bf cov}(X_i,\overline{X})\Big)\to 0.
\]

Notice also that $\frac{1}{n}{\bf
Var}\sum\limits_{i=1}^{n}\varepsilon_i^{V_i}\to\sigma^2$ as
$n\to\infty$  so (2) follows.

{\bf Step 2} Let $[t]$ be the integer part of $t$. For any fixed $m$
and for $0\leq s_1<\dots<s_m\leq1$, $k_i=[ns_i]$, we establish weak
convergence, as $n\to \infty$, of vector
$\vec{\eta}=\frac{1}{\sigma\sqrt{n}}(\widehat{\Delta}_{k_1}^0,\dots,\widehat{\Delta}_{k_m}^0)$
to vector $\vec{{Z}_F^0}=(Z_F^0(s_1),\dots, Z_F^0(s_m))$.

From  (1), (2) and from convergences $\overline{(X^0)^2}\to{\bf
Var}\xi_1$ a.s.,
$\frac{1}{n}\sum\limits_{i=1}^{k_i}X_i^0\to~GL_F^0(s_i)$ a.s.
(Goldie, 1975),
it is enough to prove $\vec{\zeta}\Longrightarrow\vec{Z_F^0}$ where
$\vec{\zeta}=~\frac{1}{\sigma\sqrt{n}}(\Delta_{k_1}^0,\dots,\Delta_{k_m}^0)$,
\[
\Delta_{k_j}^0=\sum\limits_{i=1}^{k_j}\varepsilon_i^0-\frac{GL_F^0(s_j)}{{\bf
Var}\xi_1}\sum\limits_{i=1}^{n}\varepsilon_i^0{\bf
E}X_i^0=\sum\limits_{i=1}^{k_j}\varepsilon_i^0-\frac{GL_F^0(s_j)}{{\bf
Var}\xi_1}\sum\limits_{i=1}^{n}\varepsilon_i^{V_i}{\bf E}X_i^0.
\]

{\bf Step 3} We prove weak convergence
$\vec{\zeta}\Longrightarrow\vec{Z_F^0}$ using characteristic
functions. Notice that
\[
\sum\limits_{j=1}^{m}t_j
\left(\sum\limits_{i=1}^{k_j}(\varepsilon_i^{V_i}
-\overline{\varepsilon})-\frac{GL_F^0(s_j)}{{\bf
Var}\xi_1}\sum\limits_{i=1}^{n}\varepsilon_i^{V_i}{\bf
E}X_i^0\right)
\]
\[
=\sum\limits_{i=1}^{n}\varepsilon_i^{V_i}\sum\limits_{j=1}^{m}t_j\left({\bf
I}\{i\leq k_j\}-\frac{k_j}{n} -\frac{GL_F^0(s_j)}{{\bf
Var}\xi_1}{\bf E}X_i^0\right).
\]
It is well known that the finiteness of ${\bf E}\psi_1$ implies
convergence $\frac{\psi_{n:n}}{n}\to0$ a.s. and in mean for a
sequence of i.i.d random variables $\psi_1,\dots,\psi_n,\dots$ and,
more generaly, for a stationary ergodic sequence as a consequence of
the subadditive ergodic theorem (Kingman, 1968).

Applying this fact and using H{\fontencoding{T1}\selectfont
\symbol{"AE}}lder's inequality we have ${\bf E}X_i^0=o(\sqrt{n})$
uniformly in $1\leq i\leq n$.

Let $\beta_i=\sum\limits_{j=1}^{m}t_j\left({\bf I}\{i\leq
k_j\}-\frac{k_j}{n} -\frac{GL_F^0(s_j)}{{\bf D}\xi_1}{\bf
E}X_i^0\right)$. Then  $\sum\limits_{i=1}^{n}\frac{\beta_i^2}{n} \to
C_F:=\sum\limits_{j_1=1}^{m}\sum\limits_{j_2=1}^{m}t_{j_1}t_{j_2}K_F^0(s_{j_1},s_{j_2})$
a.s. and characteristic function
$\varphi_{\vec{\zeta}}(\vec{t}\,\,)$ converges to $\exp(-C_F/2)$
a.s. Then convergence of finite-dimensional distributions follows.

{\bf Step 4.} We show that
\begin{equation}\label{1}
\mbox{the family of distributions}~\{Z_n(t), 0\leq t \leq 1\}~
\mbox{is relatively compact}.
\end{equation}

Let $S_k=\sum\limits_{i=1}^{k}\xi_{i:n},~k=1,\dots,n,~S_0=0$.

By Prokhorov's theorem (section 1 \textsection 6 in Billingsley,
1968) it suffices to show that the family of distributions of random
processes
$\left\{\frac{\widehat{\Delta}_{[nt]}^0}{\sigma\sqrt{n}},0\leq
t\leq1\right\}$, $n=1,2,\dots$, is tight. Put $k=[nt]$ and let
\[
\widehat{\Delta}_k=\sum_{i=1}^{k}\left(\varepsilon_i^{V_i}
-\frac{\overline{X^0 \varepsilon^0}}{\overline{(X^0)^2}}X_i\right).
\]
Then
$\widehat{\Delta}_k^0=\widehat{\Delta}_k-\frac{k}{n}\widehat{\Delta}_n$.
The invariance principle (e.g., part 1 of chapter 19 in Borovkov,
1998) implies tightness of the family
$\left\{\frac{\sum_{i=1}^{k}\varepsilon_i^{V_i}}{\sigma\sqrt{n}},0\leq
t\leq1\right\}$. So, to show (3), it is enough to establish
tightness of
\[
\left\{\frac{\overline{X^0
\varepsilon^0}\sqrt{n}}{\sigma\overline{(X^0)^2}}\frac{S_k}{n},0\leq
t\leq1\right\}.
\]

In turn, by Theorem 8.3 (Billingsley, 1968), it  suffices to prove
that, for any $\varepsilon>0,~\alpha>0$, there are
$0<~\delta<~1,~n_0\in{\bf   N}$ such that
\begin{equation}\label{5}
\frac{1}{\delta}{\bf   P}\left\{\sup\limits_{t\leq s\leq
t+\delta}{}\left|\frac{\overline{X^0
\varepsilon^0}\sqrt{n}}{\sigma\overline{(X^0)^2}}\frac{S_{[ns]}-S_{[nt]}}{n}\right|\geq\varepsilon\right\}\leq
\alpha,
\end{equation}
for all $n>n_0,~0\leq t\leq 1$.

Notice that $ \frac{\overline{X^0
\varepsilon^0}\sqrt{n}}{\sigma\overline{(X^0)^2}}\Longrightarrow\frac{\zeta}{\sqrt{{\bf
Var}\xi_1}}, $ and (Goldie, 1977) $ \sup\limits_{t\leq s\leq
t+\delta}{}\left|\frac{S_{[ns]}-S_{[nt]}}{n}\right|
{\to}\sup\limits_{t\leq s\leq t+\delta}{}|GL_F(s)-GL_F(t)|~{\rm
a.s.}. $ Here $\zeta$ is a standard normal random variable and
$GL_F(x)$ is the general Lorenz curve.

 By Cauchy-Bunyakowsky inequality,
\[
\sup\limits_{t\leq s\leq
t+\delta}{}|GL_F(s)-GL_F(t)|
\leq\sup\limits_{t\leq s\leq
t+\delta}{}\int_{t}^{s}|F^{-1}(x)| dx \leq\sqrt{\delta{\bf
E}\xi_1^2}.
\]

So one may choose a positive $\delta$ that satisfies (4).

{\bf Step 5.} It remains to prove
$\widehat{\sigma^2}\stackrel{\bf{p}}{\to}\sigma^2$. Indeed,
\[
\overline{\widehat{\varepsilon}^2}=\frac{1}{n}\sum\limits_{i=1}^{n}
\left(\varepsilon_i^{V_i}-\overline{\varepsilon}-\frac{\overline{X^0
\varepsilon^0}}{\overline{(X^0)^2}}(X_i-\overline{X})\right)^2
=\overline{(\varepsilon^0)^2}-\frac{(\overline{X^0
\varepsilon^0})^2}{\overline{(X^0)^2}}\stackrel{\bf{p}}{\to}\sigma^2.
\]
This completes the proof of Theorem 1.

\section*{Acknowledgements}
The work was supported in part by Russian Foundation of Basic Researches (grant 13-01-00661),
Novosibirsk State Technical University (project  2.1.1-2013),  Government of Novosibirsk Region, Russia.

Authors would like to thank Sergey Foss for many useful
discussions, Evgeny Baklanov, Alexander Sakhanenko and
Yuliana Linke for useful references.

\footnotesize
{\sc Aue, A., Horvath, L., Huskova, M., Kokoszka, P.}, 2008. Testing
for change in polynomial regression. Bernoulli 14, 637--660.

{\sc Billingsley, P.}, 1968. Convergence of Probability Measures.
New York: John Wiley $\&$ Sons.

{\sc Bischoff, W.}, 1998. A functional central limit theorem for
regression models. Ann. Stat. 26, 1398--1410.

{\sc Brown, R. L., Durbin, J., Evans, J. M.}, 1975. Techniques for
testing the constancy of regression relationships over time. J. R.
Statist. Soc. 37, 149-192.

{\sc Davydov, Y., Zitikis, R.}, 2004. Convex rearrangements of
random elements. Fields Institute Communications 44, 141--171.

{\sc Gastwirth, J. L.}, 1971. A general definition of the Lorenz
curve. Econometrica 39, 1037--1039.

{\sc Goldie, C. M.}, 1977. Convergence theorems for empirical Lorenz
curves and their inverses. Adv. Appl. Prob. 9, 765--791.

{\sc Hoeffding, W.}, 1953. On the distribution of the expected
values of the order statistics. Ann. Math. Statist. 24, 93--100.

{\sc Kingman, J. F. C.}, 1968. The ergodic theory of subadditive
stochastic processes. J. R. Statist. Soc. 30, 499--510.

{\sc Kovalevskii, A.}, 2013. A regression model for prices of
second-hand cars. Applied methods of statistical analysis.
Applications in survival analysis, reliability and quality control,
124--128.

{\sc MacNeill, I. B.}, 1978. Limit processes for sequences of
partial sums of regression residuals. Ann. Prob. 6, 695--698.

{\sc Stute, W.}, 1997. Nonparametric model checks for regression.
Ann. Statist. 25, 613--641.

\end{document}